\documentclass[12pt]{article}
 \usepackage[english]{babel}
\usepackage{color}

\usepackage[final]{epsfig}
\usepackage{amsfonts, epsfig, float}
\usepackage{graphicx}
\usepackage{amssymb,latexsym}
\usepackage{amsmath}
\usepackage{stmaryrd}
\usepackage{textcomp}
\usepackage{epsfig, float}

\usepackage{amsbsy}
\usepackage{amsthm}

\newtheorem{proposition}{Proposition}

\newtheorem{remark}{Remark}

\newtheorem{lemma}{Lemma}

\newtheorem{corollary}{Corollary}
\newcommand{\be}{\begin{equation}}
\newcommand{\ee}{\end{equation}}
\newcommand{\bea}{\begin{eqnarray}}
\newcommand{\eea}{\end{eqnarray}}
\newcommand{\bst}{\boldsymbol{\xi}}

\newcommand{\cop}{{\cal COP}(p)}

\begin{document}

\title{Representation of Zeros  of a  Copositive Matrix via Maximal Cliques of a  Graph
\thanks{This work was partially supported by    research program ‘Convergence’ (Republic Belarus), Tasks 1.3.01 and 1.3.04, and by
 Portuguese funds through CIDMA - Center for Research and Development in Mathematics and Applications, and FCT -
 Portuguese Foundation for Science and Technology, within the project  UIDB/04106/2020,  doi.org/10.54499/UIDB/04106/2020}}

\author{Kostyukova O.I.\thanks{Institute of Mathematics, National Academy of Sciences of Belarus, Surganov str. 11, 220072, Minsk,
 Belarus  ({\tt kostyukova@im.bas-net.by}).}  Tchemisova T.V.  \thanks{University of Aveiro, Campus Universitário de Santiago, 3800-198, Aveiro, Portugal  ({\tt tatiana@ua.pt})}}

\date{}

\maketitle

\begin{abstract}   There is a profound connection between copositive matrices and graph theory.
 Copositive matrices    provide a powerful tool for formulating and solving various challenging graph-related problems.
  Conversely, graph theory provides a rich set of  concepts   and techniques
that can be applied to  analyze key properties of copositive matrices, including their eigenvalues and spectra.

  In  this paper, we present  new aspects of the relationship between  copositive matrices and graph theory.
 Focusing on  the set of normalized  zeros of a copositive matrix, we investigate its properties and demonstrate that this set
  can be   expressed  as a union of   convex hulls of subsets of minimal zeros.
 We show  that  these subsets are   connected with the  set of  maximal  cliques of a special graph constructed on the basis  of   the set of
	minimal zeros  of this matrix. We develop an algorithm for constructing
	both the set   of  normalized minimal zeros and the set of all normalized zeros  of a copositive matrix.
		
\end{abstract}



\section{ Introduction}\label{Int}

 A symmetric matrix  producing  a non-negative quadratic form for all non-negative vectors, is called {\it copositive}.
 Copositive matrices play an important role in various fields such as optimization, game theory, economics, and operations research.
 Copositive matrices  define  the cone of copositive matrices (copositive cone)
which extends the concept of positive semidefinite matrices and provides a rich field for
theoretical exploration and practical applications.

 The importance of studying copositive matrices stems from their applicability in solving complex optimization problems,
particularly in Copositive   Programming (CP) and some combinatorial problems   arising in real-world scenarios.
The copositivity condition helps in formulating and solving these problems more efficiently, providing deeper insights
into structure of feasible and optimal solution sets. The diversity of copositive formulations of problems from different domains of optimization
 (continuous and discrete, deterministic and stochastic, robust optimization with uncertain objective, and others)
 is described in \cite{Bomze2012, Dur_2010}, {\it et al}. According to M. D\"{u}r \cite{Dur_2010}, CP is
 "a powerful modeling tool which interlinks the quadratic and binary worlds".
 Being formally very similar to that of Semidefinite Programming (SDP), the copositive programs are {\it NP-hard} since testing copositivity
of matrices is {\it co-NP-complete}.

 There is a  strong connection  between  the  theory of copositive matrices and graph theory.
   Numerous  problems  in graph theory  can be
 expressed  in terms of copositive matrices,   while  many properties of
  these  matrices can  be exploited  using
 results    from graph theory. Copositive matrices  are  used   for  formulating
and solving  the problem of finding the clique number of a graph
 (the number of vertices in  its largest  clique).
     This relationship further extends to the study of stability (independence) numbers, chromatic
 numbers, and other graph invariants, where copositive matrices provide
 valuable insights and methods for tackling these problems (see \cite{Dick2021, Klerk, Povh2007, Povh2013,  Vargas}, {\it et al}),
making copositive matrices  a powerful tool  of combinatorial optimization.
  Additionally,  the interplay between copositive matrices and graph
	theory enriches both fields, fostering new methodologies and applications
	in combinatorial and quadratic optimization. The results of this paper illustrate one more aspect of this relationship.

 An important direction in the study  of copositive matrices is the investigation of their {\it zeros} and {\it minimal zeros}. Given a copositive matrix $X$,
its zero   $\tau$ is a  nonzero nonnegative vector such that $\tau^\top X\tau =0.$
A zero $\tau$ of $X$ is {\it  minimal} if there does not exist another zero $t$ of $X$ such that its support, denoted here as ${\rm supp}(t)$, is a strict subset of ${\rm supp}(\tau)$
(see \cite{H1}).

	Understanding   properties   of  matrix  zeros is essential, as  the zero sets of a copositive matrices  are  used in the studying the  structure
			of the copositive cone \cite{H1, KT_Opt_2022}.
 Minimal zeros, in particular, hold special importance.  For the references see, {\it e.g.} \cite{Dick2021, Dick2016, Hild2017, H1}.
  For  conic optimization problems over  the  copositive cone,   research into the characterization and representation of zeros and minimal zeros of copositive matrices
   helps  identify the most restrictive conditions under which a feasible solution exists, thereby aiding in the development of efficient optimization algorithms.

  The purpose  of this paper is to further study the properties of normalized zeros and minimal zeros of a given copositive
 matrix and develop an algorithmic procedure for  constructing these zeros.

The paper is organized as follows.
 Section 1 provides an introduction.
 In section 2, we  introduce  the  main  notation and formulate the problem  addressed  in the paper.
In this section,  given a copositive matrix $X$,  we define the set    of all its normalized zeros and discuss the structure of this set.
In section 3,   we study the normalized minimal zeros of $X$.  It is shown that the support of a minimum zero should satisfy   certain  two conditions. Based on these conditions, an algorithm for constructing the sets  of all minimal zeros of $X$  is elaborated.
In section 4, for  the given  matrix $X$, we introduce an undirected graph $G=G(X)$,
 referred to as the {\it minimal zeros graph}, using the corresponding
 extended minimal zeros support set.
We establish a correspondence between the   structure of the   subsets defining the set   of all normalized zeros of $X$  and the set of all maximal cliques of this graph.  Using this correspondence, we  derive several properties of zeros of  the   matrix  $X$ and  the  maximal cliques  of  the corresponding graph  $G(X)$.
In particular, we show that
 the set of all normalized zeros of  $X$   can be uniquely represented as
 a union of   the convex hulls of     a finite number of   subsets  of the set of  the   minimal  zeros defined by the maximal cliques  of   G(X). In section 5,
	it is shown that   the same undirected  graph may be  generated  by  different copositive matrices through  the
	corresponding  extended minimal zeros support sets (which may be also different). A simple example   illustrates such a case.

\section{Problem statement  and notation }

 For  a given   natural number $p\geq 2$, let   $\mathbb R^p$   be
the $p$ - dimensional Euclidean vector space with  the  standard orthogonal basis $\{e_k, k=1,\dots,p\}$ and $\mathbb S(p)$    be  the space of symmetric $p\times p$ matrices.

Let $\mathbb S_+(p)$ denote the cone of symmetric  positive semidefinite $p\times p$ matrices,
$$ \mathbb S_+(p):=\{A\in \mathbb S(p):t^\top At\geq 0 \ \forall t \in \mathbb R^p\},$$
and $\cop$ denote the cone   of   symmetric  copositive  $p\times p$ matrices  which can be defined as
$$ \cop:=\{A\in \mathbb S(p):t^\top At\geq 0 \ \forall t \in \mathbb R^p,\ t\geq {\bf 0}\}$$
 or equivalently, as
$$\cop:=\{A\in \mathbb S(p):t^\top At\geq 0 \ \forall t \in T\},$$
where $T:=\{t\in  \mathbb R^p: t \geq {\bf 0}, ||t||_1=1\}$  is a simplex in $\mathbb R^p$.  Here ${\bf{0}}$ denotes  the null  $p$- vector  and $t \geq {\bf 0}$ iff each component of vector $t$ is non-negative.

 Set  $P:=\{1,\dots,p\}$. Given  a vector $t=(t_k, k \in P)^\top\in  \mathbb R^p$, $ t \geq {\bf 0}$,  denote  by
${\rm supp}(t)$  its {\it support}: ${\rm supp}(t):=\{k \in P: \ t_k>0\}$.

Suppose that a matrix $X\in \cop $ is given. Denote by $T_0$ the set of all {\it normalized zeros}   of  $X$:
$$T_0:=\{t\in T: \, t^\top Xt=0\}.$$
A zero $\tau\in T_0$  of the matrix $X$  is called   {\it minimal}  if there does not exist another   zero  of $X$,   $t\in T_0$,
 such that  and ${\rm supp}(t)\subset {\rm supp}(\tau).$  It is known that for $X\in \cop$, there is a finite number of normalized minimal zeros.  In \cite{KT_Opt_2022}, it is shown that the minimal zeros of certain matrices play an important role in explicit representation of the faces of the copositive cone $\cop $.

It is a known  fact (see {\it e.g.} \cite{Baumert66, Dick2016}),  that if $\tau\in T_0$,  then $X\tau\geq {\bf 0}$, wherefrom we conclude that
\be t^\top X\tau\geq 0\ \forall \tau\in T_{0},\ \forall  t \in T.\label{1Z}\ee

In \cite{  KTD_2020} (see Lemma 1, where  a convex subset ${\cal A}(x)\subset \cop$ is replaced by $X\in \cop$), it is shown   that the set $T_0$ is  either  empty
 or  can be presented as a union  of a finite number polytopes. Hence, in the case
 $T_0\neq \emptyset$, we have
\be T_0:=\bigcup\limits_{s\in \widetilde S}\widetilde T_0(s) \mbox{, where }  \ \widetilde T_0(s),\ s \in \widetilde S, \mbox{ are
     polytopes}.\label{z1}\ee
 Here $\widetilde S$ is a finite  subset of $\mathbb N$.
Such a  representation is  useful  because, although the set $T_0$  typically  consists  of an infinite number of elements,
   formula (\ref{z1})  allows one to   describe $T_0$ using a finite set of data,
 namely,  the  set of   vertices of  the   polytopes  $\widetilde T_0(s)$, $ s \in \widetilde S.$

In general,   representation  (\ref{z1}) is not unique. In what follows, we are interested  in
 representations having   the minimal number  of  polytopes.   A representation  (\ref{z1})
with the smallest value of  $|\widetilde S|$,   will be called  the  {\it minimal} representation of the set $T_0 $.
 For the minimal representation, a finite data set allowing  one  to describe the set
$T_0$,  is minimal (by number of elements).

  The  aim of the paper  is   to    light on  new   properties of normalized zeros of  a given matrix $X\in \cop$
	and describe  the  rules for constructing   the   minimal  representations of the corresponding   set $T_0$.

\section{Constructing   the set of normalized  minimal zeros}

Consider a copositive matrix  $X\in \cop $  and the corresponding set $T_0$ of all its normalized zeros. Suppose that  $T_0\not=\emptyset.$  Let us denote by
 $$\{\tau(j), j \in J\} , \ J\subset \mathbb N,\ |J|<\infty,$$
the set of all normalized minimal zeros  $\tau(j), j \in J,$ of $X.$

In what follows, we will show  that  there exists a finite set $ \{J(s),s \in S\}$ of subsets of the index set $J$  such that
 the following representation of the set $T_0$ takes place:
$$ T_0:=\bigcup\limits_{s\in S}T_0(s),  \mbox{ where }  \ T_0(s)={\rm conv} \{\tau(j), j \in J(s)\}, \ s \in S.$$
 Thus, it is  clear  that  the  minimal zeros  $\tau(j), j \in J,$ of $X$   play an important role in
  describing of the set $T_0.$ In this  section,  we  will  study  the  properties of these  zeros, and
 use  these  properties    to describe   rules   for their  construction.

 Let $\overline P$ be  a subset  of the set $P=\{1,\dots, p\}$:  $\overline P\subset P$.
 Denote by $X(\overline P)$ the  corresponding principal sub-matrix of    the  $p\times p $\emph{} matrix $X$:
\be\label{fX} X(\overline P)=\begin{pmatrix}
X_{kq}, q\in \overline P\cr
k\in \overline P\end{pmatrix},\ee
where $X_{kq}$  {denotes the  element } of $X$  on the position $(k, q)$.

 The following lemma  and proposition  describe some properties of  the support of a minimal zero   of a copositive matrix
through   properties of     its   principle  sub-matrix.

\begin{lemma} \label{L1} (\cite{H1}, Lemma 3.7)  Given  $X \in \cop$, let  $\overline P\neq \emptyset$  be a subset of the set $P$. Then $\overline P$
is the support of a minimal zero of $X$ if and only if the  principal  sub-matrix $X(\overline P)$  is positive
semidefinite with corank 1, and  its null space (kernel) is generated by a  vector with  all  positive elements.
\end{lemma}

 The next proposition is  a slight modification of Lemma \ref{L1}.  The conditions formulated in this proposition are better suited for software implementation.

\begin{proposition}\label{P1}
 Let $X \in \cop$ and  $ P=\{1,\dots, p\}.$

 I. Suppose that the set $\ \emptyset \neq\overline P\subset P$ satisfies the following conditions:
\begin{enumerate}
\item[(A)] $\ \ \ {\rm rank}(X(\overline P))=|\overline P|-1;$
\item[(B)] $\ \ \ \exists \, \bar t =(\bar t_k,k \in \overline P)$ such that
$X(\overline P)\bar t ={\bf 0},\ \bar t>{\bf 0}, \ ||\bar t ||_1=1,$
 where the submatrix $X(\overline P)$ is constructed by formula (\ref{fX}).
\end{enumerate}
Then the vector
\be \bar{\tau}=(\bar{\tau}_k =\bar t_k ,\, k \in \overline P,\ \bar{\tau}_k(j)=0, \, k \in P\setminus \overline P) \label{tauj}\ee
is a normalized minimal zero of $X$.

 II.  Let $\bar{\tau}=(\bar{\tau}_k, k \in P)$ be a minimal zero of $X$. Then the set $\overline P:={\rm supp}(\bar{\tau})$ satisfies
conditions  {$(A)$} and  {$(B)$}.

\end{proposition}

{\bf Proof.}
 First, let us prove  statement I of the proposition.  Suppose that   conditions $(A)$ and $(B)$ hold true. Let us show that the vector
 $\bar{\tau}$  constructed  {as in } (\ref{tauj}) is a normalized minimal zero  of   matrix  $X$. First of all,  note that

$\bullet$
  conditions $X\in \cop$ and  (B) imply
\be  X(\overline P)\in \mathbb S_+(|\overline P|),\label{1.1}\ee

$\bullet$ it follows from  conditions $(A)$ and $(B)$ that  the  null space  of $X(\overline P)$ is generated by the vector $\bar t>{\bf 0}$.

Let  $\bar{\tau}$ be a vector constructed by the rule (\ref{tauj}).
 Then, taking into account the properties mentioned above and Lemma \ref{L1}, we conclude that  $\bar{\tau}$
is a minimal normalized  zero of $X$ and  statement I is proved.

To prove statement II of the proposition,  let    $\bar{\tau}$  be a minimal normalized  zero of $X$. Set
 $\overline P:={\rm supp}(\bar{\tau})$ and
denote $\bar t:=(\bar t_k =\bar{\tau}_k, k \in \overline P)$. By construction, the vector $\bar t$ satisfies condition (B).
 Then   condition $(A)$ follows from Lemma \ref{L1}.  $\ \Box$

\begin{proposition}\label{P10}  Let  $X\in \cop$ and a  subset $\overline P \subset P$ satisfy  conditions $(A)$ and $(B)$.  Then
$${\rm det}X( \overline P \setminus \{i\})\not=0\ \forall i \in \overline P.$$
\end{proposition}
 {\bf Proof.}  Suppose the contrary:  for some  $i_0\in  \overline P$, it holds   ${\rm det}X(\overline P\setminus \{i_0\})=0$. Then  there exists
$\mu=(\mu_k, k\in \overline P\setminus \{i_0\})^\top\not={\bf 0}$ such that $X(\overline P\setminus \{i_0\})\mu={\bf 0}$. This implies
$\bar \mu^\top X(\overline P)\bar \mu=0\ \mbox{ for } \bar\mu:=(\bar \mu_k=\mu_k, \, k \in \overline P\setminus \{i_0\}, \bar \mu_{i_0}=0)^\top,$
wherefrom, taking into account  condition (\ref{1.1}), we obtain  $ X(\overline P)\bar \mu={\bf 0}.$

 It follows from  condition  (B) that there exists $ \bar t =(\bar t_k,k \in \overline P)$ such that
$X(\overline P)\bar t ={\bf 0},\ \bar t>{\bf 0}.$ Since $\bar t_{i_0}>0$, $\bar \mu_{i_0}=0$, and $\bar \mu\not={\bf 0}$, it is evident that the vectors
$\bar t$ and $\bar \mu$ are linearly independent.
Thus we have  obtained that  for the pair of linear
independent vectors $\bar t$ and $\bar \mu$, it holds  $X(\overline P)\bar t ={\bf 0},$ $X(\overline P)\bar \mu ={\bf 0}$, which  contradicts    condition $(A)$. $\ \Box$

\vspace{2mm}

 Recall that the  normalized   minimal zeros of the matrix $X$ are   indexed by elements
from a finite set  $J$, and  denoted as  $\tau(j), \ j\in J.$  Based on  Proposition \ref{P1}, we  conclude
that, in order to determine   the set of  all minimal zeros  of $X$,  we  have to find  in $P$  all subsets
$P(j), j\in J,$   satisfying   conditions $(A)$ and $(B)$.

\begin{corollary}\label{WWW} Let ${\{}P(j),  \  j \in J{\}}$ be the set of all different subsets  of the set
$P$  satisfying    conditions $(A)$ and $(B)$. Then  the following conditions are   satisfied:
\be P(i)\not \subset P(j)\ \forall i \in J, \ \forall j \in J, \ i\not=j.\label{6}\ee
\end{corollary}
{\bf Proof.} Suppose that, on the contrary, there exists a pair of indices $i_0 \in J$, and
  $ j_0 \in J,$   $ i_0\not=j_0,$ such that  the inclusion
$P(i_0)\subset P(j_0)$   takes place.  Notice that  all sets $P(j), $ $ j \in J,$
are different, and hence  $P(i_0)\not =P(j_0).$

 By assumption, the sets $P(i_0)$ and $ P(j_0)$ satisfy  the  conditions $(A)$ and $(B)$. From
Proposition \ref{P1}, it   follows  that  for  the set $P(i_0)$, there exists  a minimal zero $\tau(i_0)$ such that ${\rm supp}(\tau(i_0))=P(i_0),$
and  for  the set $P(j_0)$, there exists  a minimal zero $\tau(j_0)$ such that ${\rm supp}(\tau(j_0))=P(j_0).$  Hence,  for two different
minimal zeros $\tau(i_0)$ and $\tau(j_0)$,  we obtain
 ${\rm supp}(\tau(i_0))=P(i_0)\subset P(j_0)={\rm supp}(\tau(j_0)).$ But the inclusion
${\rm supp}(\tau(i_0))\subset {\rm supp}(\tau(j_0))$ and  the  inequality $\tau(i_0)\not=\tau(j_0)$ contradict the statement that $\tau(i_0)$ is a minimal zero
of $X$. $\ \Box$


\vspace{2mm}

 For a given subset $\overline P\subset P$,  consider the corresponding   vector ${\bf a}(\overline P)=(a_k(\overline P), k \in P)^\top, $ defined by the rules
 $$ a_k(\overline P)=1\mbox{ if } k\in   \overline P,  \ a_k(\overline P)=0\mbox{ if } k\in P\setminus   \overline P.$$

Then condition (\ref{6}) is equivalent to the following one:
$$ ( {\bf a}(P(i)))^{\top} {\bf a}(P(j))<\min\{|P(i)|, |P(j)|\} \  \forall i \in J, \ \forall j \in J, \ i\not=j.$$

 At the end of this section, we will describe an algorithm for constructing the sets $P(j)\subset P,$ $ j \in J,$  that satisfy conditions $(A)$ and $(B)$.

\vspace{6mm}

\centerline{{\bf Algorithm.}}

\vspace{2mm}

{\bf Step 1.}   Set  $\Delta I_1:=\{k\in P:X_{kk}=0\}.$

If $\Delta I_1\not =\emptyset,$ then set:
$P(1,i){:}=e_i, \  {\bf a}(1,i){:}=e_i, \ i\in \Delta I_1,$  where $\{e_k, k=1,\dots,p\}$
 is the standard orthogonal basis the space $\mathbb R^p$. Replace $P$ by  $P\setminus \Delta I_1$ and go to Step 2  with  the  new set  ${P}$.

If $\Delta I_1 =\emptyset,$ then  the set $P$ does not change.

\vspace{2mm}

{\bf Step 2.}   Let  $P(2,i),$ $ i \in I_2,$   with  $I_2=\{1,\dots,C^2_p\},$
$C^2_p=\frac{p!}{2!(p-2)!}$,  be    all different subsets of   $P$ consisting of two elements.
Denote:
$$\Delta I_2:=\{i \in I_2: P(2,i) \mbox{ satisfies  conditions (A) and (B)}\},$$
and construct  the vector  $ {\bf a}(2,i):={\bf a}(P(2,i))$ corresponding to $P(2,i)$ for $i \in \Delta I_2.$

Go to Step 3.

\vspace{2mm}

{\bf Step $k$, $ k\geq 3.$} At  the beginning of this step,  we have    vectors
$${\bf a}(m,i), \,i \in \Delta I_m,\ m=2,\dots,k-1.$$
 Note that it  can occur that   $\Delta I_m=\emptyset$ for some $m\in \{2,\dots,k-1\}.$

 Let $P(k,i),$ $ i \in I_k,$  with $I_k=\{1,\dots,C^k_p \},$  $C^k_p=\frac{p!}{k!(p-k)!}$, be all different subsets of $P$ consisting of $k$ elements. Denote
$ {\bf a}(k,i):= {\bf a}(P(k,i))$,  $ i \in I_k, $ and set
\begin{equation*}\begin{split}&\bar I_k:=\{i \in I_k: ({\bf a}(k,i))^{\top}{\bf a}(m,j)<m\ \forall j\in \Delta I_m, \ m=2,\dots,k\},\\
&\Delta I_k:=\{i \in \bar I_k: P(k,i) \mbox{ satisfies (A) and (B)}\}.\end{split}\end{equation*}

 If $k<p$, go to Step $k+1$. If $k=p,$ then STOP.

\vspace{3mm}

Note that if    condition $(A)$ holds true for $\overline P\subset P$, then  condition $(B)$ can be tested by the following rules.

Choose any index $i_*\in \overline P.$ If ${\rm det }X(\overline P\setminus i_*)=0$, then it follows from Proposition \ref{P10}
that    condition $(B)$ can not be fulfilled.

Suppose  that
${\rm det }X(\overline P\setminus i_*)\not=0.$ Calculate the vector
$$y(\overline P, i_*)= (y_k(\overline P,i_*), k \in \overline P\setminus i_* )^\top =-(X(\overline P\setminus i_*))^{-1}(X_{k\,i_*}, k \in \overline P\setminus i_*)^\top.$$
If $y(\overline P, i_*)> {\bf 0}$, then   condition $(B)$ is satisfied, otherwise  no.

\vspace{4mm}

Note that if   conditions $(A)$ and $(B)$ are fulfilled for  $\overline P\subset P$, then the vector $ {\bst}(\overline P)=( {\bst}_k(\overline P), \, k \in P)$, constructed by the rule
\be  {\bst}_k(\overline P)=\left\{\begin{array}{ll}
y_k(\overline P, i_*)/\bar \beta &\mbox{ if } k \in \overline P\setminus \{i_*\},\cr
 1/\bar \beta &\mbox{ if } k=i_*,\cr
0&\mbox{ if } k \in P\setminus (\bar  P\cup \{i_*\}),
\end{array}\right. \mbox{ with } \bar \beta=1+||y(\overline P, i_*)||_1\label{rule},\ee
is a  normalized  minimal zero of $X$ corresponding to $\overline P$.

\vspace{3mm}

Let $ {\bf a}(k,i)=(a_q(k,i),q \in P),$ $ i \in \Delta I_k,\,k=1,\dots,p,$ be  the vectors  constructed  at the iterations of   the algorithm. Then the
set $\{P(k,i),\  i \in \Delta I_k,\,k=1,\dots,p\}$, whose
 elements are the sets $P(k,i){:}=\{ q \in P: \   a_q(k,i)=1\} $    for all $i \in \Delta I_k,\,k=1,\dots,p$,  gives us
 the  sought-after  set $\{P(j), j \in J\},$ $|J|=\sum\limits_{k \in P}|\Delta I_k|,$ of  all subsets of $P$
 satisfying     the conditions  $(A)$ and $(B)$. Having the sets $P(j), j \in J,$
 it is  straightforward to find the corresponding set $\{\tau(j), j \in J\} $ of all  normalized  minimal zeros of $X$
  using the rules  {in} (\ref{rule}):
$\tau(j)= {\bst}(P(j)), $ $j\in J.$

\begin{remark}   It is worth noting  that   by utilizing certain  properties of minimal zeros (see, for example, Theorem 3.11 and Corollary 3.12 \cite{H1}),
 one can  potentially  further  reduce the search for subsets in the algorithm described above.
\end{remark}
\begin{remark}
 In  \cite{Dic2019},  a   novel  method   for  certifying  the copositivity of  any  given  matrix
 was  introduced.   This certificate  is constructed
by solving a finite  number of  linear systems and can be  subsequently verified by checking a finite number of  linear inequalities. Additionally, this certificate can  be
  employed to generate the set of minimal zeros of a copositive matrix.
\end{remark}

\section{  Construction of the minimal representation  for  the set of all normalized zeros}

In this section,   for a given matrix $X\in \cop$,  we  will study the properties   of    the corresponding   set $T_0$  of all its
normalized zeros   using the results from the graph theory.

Suppose that for a given matrix $X\in \cop,$  all corresponding normalized  minimal zeros $\tau(j),$ $ j \in J,$ are found. Let us denote
\be    M(j):=\{k \in P:e^\top_k X\tau(j)=0\},\ j \in J.\label{N1}\ee

 In  \cite{Hild2020}, for a given  $X\in\cop$, the set $ \ \mathcal{E}=\{({\rm supp}(\tau(j)),M(j)), j \in J\}$,  consisting  of   pairs  of subsets of $P$ constructed
  using the corresponding set of normalized minimal zeros $\{\tau(j), j \in J\}$, is   referred to as  {\it the extended minimal zeros support set} of $X.$

Using the   extended minimal zeros support set of $X,$ define
  the set $V$ of   ordered pairs of indices   from  the set $J$  as follows:
  \be V=\{ (i,j):i\in J,\, j\in J, \, i<j,\,{\rm supp}(\tau(i))\subset M(j)\}\label{Ext}\ee
 and  consider  the corresponding  undirected  graph $G(X)=(J, V)$
 with    the  set of vertices    $J$   and   the set of  edges $V$.
Notice that  the  graph $G(X)$ is uniquely defined by the extended minimal zeros support set of $X.$

 In what follows, we will refer to this graph as to  the   {\it minimal zeros graph} of matrix $X$.

\begin{remark}
 It is   straightforward   to see that the set $V$ of  ordered pairs of indices  from   the set $J$ defined in (\ref{Ext})
can  alternatively  be expressed  as follows:
 \be V= \{(i,j): \ i \in J, \ j \in J,\ i< j,\ (\tau(i))^\top X\tau(j)=0\}.\label{N2}\ee
\end{remark}

 A {\it clique}  $I$ of a  graph $G$ is a subset of the vertex set, $ J$,
such that every two distinct vertices in $I$   are adjacent.
A {\it maximal clique} of $G$  is a clique that is not a proper subset of any other clique of $G$.
 In what follows, for $j\in J,$ we consider that   the  singleton  set  $\{j\}$  is  (trivially)  a clique  of  $G$.

Let \be \{J(s), s \in S\}\label{cc1}\ee be the set of all (distinct) maximal cliques  of   the graph  $G(X)$.
It is  evident that if $\overline J\subset J$ is a  clique of  $G(X)$, then there exists $\bar s \in S$ such that $\overline J\subset J(\bar s).$

\begin{proposition}\label{PP2} Given a matrix $X \in \cop$, consider the    set (\ref{cc1})  of all maximal cliques of
 the corresponding minimal zeros  graph $G(X)$. Then
 the set $T_0$  of all normalized zeros of $X$ admits the representation
\be T_0=\bigcup\limits_{s \in S}T_0(s), \mbox{ where } T_0(s):={\rm conv}\{\tau(j), j \in J(s)\},\ s \in S.\label{Ta}\ee
  Here $\{\tau(j), j \in J\}$ is the set of all normalized zeros of $X.$
\end{proposition}
{\bf Proof.}   Let  $t\in T_0(s)$  for  some  $s \in S$.   Then   $t$ admits a representation
$$t=\sum\limits_{j \in J(s)}\alpha_j\tau(j),\ \alpha_j\geq 0,\, j \in J(s), \ \sum\limits_{j \in J(s)}\alpha_j=1.$$
Since  $J(s)$ is a clique  of  the minimal zeros graph  $G(X)$, we have
\be (\tau(i))^\top X\tau(j)=0, \ i \in J(s),\ j \in J(s),\label{10.1}\ee
 {wherefrom, taking into account the representation of the vector $t$, we get}
$$t^\top Xt=\Big(\sum\limits_{i \in J(s)}\alpha_i\tau(i)\Big)^\top X\Big(\sum\limits_{j \in J(s)}\alpha_j\tau(j)\Big)=
\sum\limits_{i \in J(s)}\sum\limits_{j \in J(s)}\alpha_i\alpha_j(\tau(i))^\top X\tau(j)=0.$$
 Then, by definition,    $ t \in T_0,$ and consequently, $T_0(s)\subset T_0$. Hence
\be \bigcup\limits_{s \in S}T_0(s)\subset T_0.\label{10}\ee

Now, let us consider any $t\in T_0.$ From Corollary 3.4 in \cite{H1}, it follows that there exists a set $J_*\subset J$ such that
\be t=\sum\limits_{j\in J_*}{\alpha_j}\tau(j) \mbox{ with some } \alpha_j> 0, \, j \in J_*. \label{11}\ee

Since $t \in T_0,$ then $t^\top Xt=0$ and hence,
\be t^\top Xt=\Big(\sum\limits_{i \in J_*}\alpha_i\tau(i)\Big)^\top X\Big(\sum\limits_{j \in J_*}\alpha_j\tau(j)\Big)=
\sum\limits_{i \in J_*}\sum\limits_{j \in J_*}\alpha_i\alpha_j(\tau(i))^\top X\tau(j)=0.\label{2Z}\ee
Notice that due to (\ref{1Z}), we have
$$(\tau(i))^\top X\tau(j)\geq 0 \ \forall i \in J,\ \forall   j\in J.$$
From the inequalities above, the inequalities $\alpha_i\alpha_j>0,$ $i \in J_*,$ $j \in J_*,$ and
equality (\ref{2Z}), one can conclude that for all $ i \in J_*,\   j \in J_*,$ such that $\ i<j$, it holds   $(\tau(i))^\top X\tau(j)=0 $, and hence, $(i,j)\in V.$
Consequently, $J_*$ is a clique  of $G(X)$  and there exists $\bar s\in S$ such that $J_*\subset J(\bar s).$ Then it follows from (\ref{11}) that
$$t \in T_0(\bar s)\ \Longrightarrow\ t \in \bigcup\limits_{s \in S}  T_0( s)\ \Longrightarrow\ T_0\subset \bigcup\limits_{s \in S}  T_0( s).$$
The latter inclusion and (\ref{10}) imply that $T_0= \bigcup\limits_{s \in S}  T_0( s).$ $\ \Box$

Based on  Proposition \ref{PP2},  we can deduce the following   useful properties    of a copositive matrix $X$,  its zeros,   and the corresponding minimal zeros graph $G(X).$

\begin{corollary}\label{C00}  The set $\{\tau(j), j \in J\}$ of all normalized zeros of $X \in \cop$ coincides with the set of all
 vertices of  the set  $\ {\rm conv}\, T_0.$
\end{corollary}
{\bf Proof.} It follows from (\ref{Ta}) that the set ${\rm conv}T_0$ is a  polytope. Denote by $\{\xi(i), i \in I\}$ the set  of all  its vertices.
 Then it follows from  the  representation (\ref{Ta}) that $\{\xi(i), i \in I\}\subset \{\tau(j), j \in J\}$. Hence, without loss of generality, we can consider that
 $\{\xi(i), i \in I\}= \{\tau(j), j \in \bar J\}$  for some $\bar J\subset J$, $|I|=|\bar J|.$

   Suppose that there exists $j_0\in J$ such that $\tau(j_0)\not \in \{\tau(j), j \in \bar J\}.$
 Since $\tau(j_0)\in T_0\subset {\rm conv}T_0,$ there exist a subset $J_*\subset \bar J$ and numbers $\alpha_j>0$,  $ j \in J_*,$
 such that $\sum\limits_{j\in J_*}\alpha_j=1$  and  $\tau(j_0)=\sum\limits_{j\in J_*}\alpha_j\tau(j).$ This implies that
 ${\rm supp}(\tau(j))\subset {\rm supp}(\tau(j_0))$ for all $j \in J_*.$ Since $\tau(j_0)$ is a minimal zero, it follows from the latter inclusions
(see Lemma 3.5  in \cite{H1})   that $\tau(j_0)=\tau(j), j \in J_*\subset \bar J.$ But this contradicts the assumption that $\tau(j_0)\not \in \{\tau(j), j \in \bar J\}.$  $\ \Box$

\vspace{2mm}

 Denote
\be P_*(s):=\bigcup\limits_{j\in J(s)}{\rm supp}(\tau(j)),\ s \in S.\label{c2}\ee

 \begin{corollary}\label{C01}    Let  $T_0$ be the set of all normalized zeros of a   matrix   $X\in \cop$. A zero $t\in T_0$ belongs to $T_0(s)$ with some
 $s \in S$ iff $\ {\rm supp}(t)\subset P_*(s).$
 \end{corollary}
 {\bf Proof.}
 Suppose that $\bar t\in T_0(s)$ for some $s\in S.$ Then $\bar t$ admits a representation
 $\bar t=\sum\limits_{j \in J(s)}\alpha_j\tau(j)$ with  certain  $\alpha_j\geq 0$ $\forall j \in J(s).$ This implies that
 ${\rm supp}(\bar t)\subset P_*(s)$.

Now let us prove that the conditions  $\bar t\in T_0$ and  ${\rm supp}(\bar t)\subset P_*(s)$ imply the inclusion $\bar t \in T_0(s).$

First, let us show that for  any $s\in S,$  the   matrix $X(P_*(s))$ is positive semidefinite.
  In fact,  consider the vector $\tau=\sum_{j \in J(s)}\tau(j)/|J(s)|.$
By construction, $\tau\in T_0(s)$ and ${\rm supp}(\tau)=P_*(s).$
It follows from these conditions and Lemma 2.4 from \cite{Dic2013} that the matrix $X(P_*(s))$ is positive semidefinite.

 \vspace{2mm}

Let $\bar t\in T_0$ be such that ${\rm supp}(\bar t)\subset P_*(s)$ with some $s\in S.$
Consider a  polytope
$$L(s):=\{t \in T:{\rm supp}(t)\subset  P_*(s),\ X(P_*(s))(t_k, \, k \in P_*(s))^\top={\bf 0}\}.$$
Let us show that
\be t\in  T_0 \mbox{ and } {\rm supp}(t)\subset P_*(s) \ \Longrightarrow\ t \in L(s).\label{Ls}\ee
In fact, since $t\in  T_0 $   and ${\rm supp}(t)\subset P_*(s)$, we have
$$0=t^\top Xt=(t_k, k \in P_*(s))X(P_*(s))(t_k, k \in P_*(s))^\top.$$  Taking into account this
 equality and  the fact that the matrix $X(P_*(s))$ is positive semidefinite,
we conclude that $X(P_*(s))(t_k, k \in P_*(s))^\top={\bf 0}.$ Thus, relations (\ref{Ls}) are proved.

 It is clear that $L(s)\subset T_0.$
Denote by $\{\mu(i), i\in N\}$ the set of all vertices of the  polytope  $L(s).$  Let us show that for $i \in N$, a normalized zero  $\mu(i)$ is a  minimal zero of
$X$.   Indeed, if suppose the contrary, then there exists a zero $\tau\in T_0$ such that the inclusion  ${\rm supp}(\tau)\subset {\rm supp}(\mu(i))$ holds
true strictly.
 It follows from this inclusion that $\tau^\top X\mu(i)=0$ and there exists $\theta\in (0,1)$ such
$\bar \tau:= \frac{1}{1-\theta}(\mu(i)-\theta \tau)\geq {\bf 0}.$
 It is easy to see  that  for the vector $\bar \tau$, it holds: $\bar \tau\in T,$ ${\rm supp}(\bar \tau)\subset {\rm supp}(\mu(i))\subset P_*(s),$ $\bar \tau^\top X\bar\tau=0.$
It follows from these relations and (\ref{Ls}) that $\tau \in L(s)$ and $\bar \tau\in L(s).$ Thus, we obtained that $\mu(i)=(1-\theta) \bar\tau+\theta \tau$, where
$\theta\in (0,1),$ $ \tau\in L(s), $  $ \bar \tau\in L(s), $ $\tau\not=\mu(i).$ But this contradicts the assumption that $\mu(i)$ is a vertex of the set $L(s)$.
Consequently, we have proved    the inclusion $\{\mu(i), i\in N\}\subset\{\tau(j), j \in J\}.$
This implies  that  $\{\mu(i), i\in N\}=\{\tau(j), j \in \bar J\}$   for  some $\bar J\subset J,$ $|N|=|\bar J|.$

  From  (\ref{Ls}), we conclude  that $\tau(j)\in L(s)$ for all $j \in J(s).$    Then, evidently,  $J(s)\subset \bar J.$

Moreover, it is easy to see that $(\tau(i))^\top X\tau(j)=0$ $\forall i \in \bar J,$   $\forall j \in \bar J,$ and hence $\bar J$ is a clique of graph $G(X)$. Taking into account that by construction, $J(s)$ is a maximal clique and  $J(s)\subset \bar J,$ we conclude that  $J(s)=\bar J,$ and consequently $T_0(s)=L(s).$

As by assumption, $\bar t\in T_0$ and  ${\rm supp}(\bar t)\subset P_*(s)$, then it follows from (\ref{Ls}) that $\bar t\in L(s)$ and, due to the equality
$T_0(s)=L(s)$,  it holds  $\bar t\in T_0(s).$ $\ \Box$

 \begin{corollary}\label{C1}    Given $X\in \cop$, let  the set  (\ref{cc1}) be the set of all maximal cliques  of  $ G(X).$ If $|S|\geq 2,$ then
\be P_*(s)\not \subset P_*(\bar s)\ \forall\, s\in S, \ \forall \,\bar s\in S,\ s\not=\bar s.\label{8.1}\ee
\end{corollary}
{\bf Proof.}  Suppose that, on the contrary, there   exist  $ s\in S, \ \bar s\in S,$   such that  $s\not=\bar s$ and  $P_*(s)\subset P_*(\bar s).$
Since $J(s)$ and  $J(\bar s)$  are different maximal cliques of $G(X)$,   it follows  that $J(s)\not\subset J(\bar s)$. Hence, there exists $i_0\in  J(s)\setminus J(\bar s)$. Then it follows from the inclusion  $P_*(s)\subset P_*(\bar s)$ that
${\rm supp}(\tau(i_0))\subset P_*(s) \subset P_*(\bar s)\subset M(j)\ \forall \, j \in J(\bar s).$
This implies that for all $j \in J(\bar s)$,   it holds
$$\tau(i_0)^\top X\tau(j)=\sum\limits_{k\in {\rm supp}(\tau(i_0))}\tau_k(i_0)e^\top_k X\tau(j)=0.$$ Hence,  from the  definition  of the set $V$ (see (\ref{N2})),  it follows  that
$ (i_0,j)\in { V}$  for all  $ j \in J(\bar s).$ 
Taking into account these inclusions   and  the fact   that $J(\bar s)$ is a clique of  the  graph $G(X)$,   one  concludes that
  the  set $\overline J(\bar s):=J(\bar s)\cup \{i_0\},$ $\overline J(\bar s)\subset J,$ $i_0\not \in  J(\bar s),$  is  a clique   of
$G(X)$ as well, which  contradicts the assumption that
$J(\bar s)$ is a maximal clique in  $G(X)$. $ \ \Box$

\begin{corollary}\label{C3} Let $\bar \tau\in T_0$ and ${\rm supp}(\bar\tau)=P_*(\bar s)$  for  some $\bar s \in S.$ Then $\bar \tau\in T_0(\bar s)$ and
$ \bar \tau\not \in T_0( s)$ for all $s \in S\setminus  \{\bar s \}.$
\end{corollary}
{\bf Proof.} Note that for any $s \in S,$
\be \mbox{the condition } \tau \in T_0(s) \mbox{ implies  the inclusion  } {\rm supp}(\tau)\subset P_*(s).\label{*11}\ee

Let $\bar \tau\in T_0$ and ${\rm supp}(\bar\tau)=P_*(\bar s)$  for  some $\bar s \in S.$   If suppose
that there exists $\hat s\in S\setminus  \{\bar s \}$ such that  $\bar \tau \in T_0(\hat s)$, then  from (\ref{*11}), one can conclude  that
${\rm supp}(\bar\tau)=P_*(\bar s)\subset P_*(\hat s)$, which  contradicts (\ref{8.1}). Hence we have shown that
$\bar\tau\not\in T_0(s)$ for all $s \in S\setminus \{\bar s \}.$ By assumption, $\bar \tau\in T_0$. Then  from    Proposition \ref{PP2}, it follows  that
$\bar\tau \in \bigcup\limits_{s \in S}T_0(s). $ This implies that $\bar \tau\in T_0(\bar s).$ $\ \Box$

\begin{corollary}\label{C2} Let  $X\in \cop$.  Consider  the set  $\{\tau(j), j \in J\}$ 
of all normalized minimal zeros of $X$  and    the set  (\ref{cc1})  of all maximal cliques   of  the graph  $G(X).$
  Suppose that $s\in S$  and  that  $\overline J$ is a clique   of  $G(X)$ such that
\be P_*(s)\subset \overline P_*:=\bigcup\limits_{j\in \overline J}{\rm supp}(\tau(j)).\label{*1}\ee
Then  the following statements hold:
$$\textbf{1})\  P_*(s)= \overline P_*,\ \   \textbf{2})\  \overline J\subset J(s),\ \ \textbf{3}) \
 \overline J\not \subset J(\bar s) \ \forall\bar s\in S\setminus  \{  s \}.$$
\end{corollary}
{\bf Proof.} Let $ \widehat J$ be a maximal clique  of  $G(X)$ such that $\overline J\subset  \widehat J.$
 Then there is $\widehat  s \in S$ such that $\widehat J=J(\widehat  s)$   and  from (\ref{*1}), it follows:
\be  P_*(s)\subset \overline P_*\subset P_*(\widehat s):=\bigcup\limits_{j\in  J(\widehat s)} {\rm supp}(\tau(j)),\label{*2}\ee
where $s \in S$ and $\widehat s\in S.$ Taking into account (\ref{8.1}), we conclude that relations (\ref{*2})  can only  hold {true}
 under the condition that $s=\widehat  s.$   In this case,  the  inclusions  given in  (\ref{*2}) imply that $P_*(s)=\overline P_*.$
 Moreover,   when  $ s=\widehat s$, we have
$\overline J\subset  \widehat J=J(\widehat s)=J(s).$

 Finally, let us show that $\overline J\not \subset J(\bar s)$ $\forall \bar s\in S\setminus   \{ s \}.$ Suppose  that on
 the contrary, there exists $\bar s\in S\setminus s$
such that $\overline J \subset J(\bar s).$ Then  $\overline P_*\subset P_*(\bar s)$ and hence,
$P_*(s)=\overline P_*\subset P_*(\bar s)$ with $s\in S,$ $\bar s \in S,$ $s\not=\bar s$, which  contradicts (\ref{8.1}).
$\ \Box$

\vspace{2mm}

 The next  proposition   shows that the  representation (\ref{Ta}) of the set $T_0$, where the
set (\ref{cc1})   consists  of all maximal cliques  of  the graph $G(X)$, is {\it a unique } minimal  representation of the set $T_0$.

\begin{proposition}\label{PP3}  Given a matrix $X\in \cop$, let   $\{\tau(j), j \in J\}$ be the set of
its normalized minimal zeros and $\{ J(s), \ s \in S\}$ be the set of all maximal cliques  of the minimal zeros graph $G(X)$.
   Suppose that the set $T_0$ admits a representation  (\ref{z1}):
$$ T_0=\bigcup\limits_{s \in \widetilde S}\widetilde T_0(s),\ \mbox{where }\widetilde T_0(s),\  s \in \widetilde S, \mbox{ are polytopes}. $$
Then   we have  $|\widetilde S|\geq |S|,$ and  the  equality $|\widetilde S|= |S|$
implies that, up to renumbering, the following holds for all $s \in S$:
 $T_0(s):={\rm conv}\{\tau(j),j \in J(s)\}=\widetilde T_0(s)$.
\end{proposition}

{\bf Proof.}  Suppose that the set $T_0$ admits a  representation   in the form  (\ref{z1}).
For $s \in \widetilde S,$ let $\{\mu(s,j),\ j \in \widetilde J(s)\}$ be the set of  all   vertices of the  polytope  $\widetilde T_0(s).$ Then
$\widetilde T_0(s)$ admits  a representation $\widetilde T_0(s)={\rm conv} \{\mu(s,j),\ j \in \widetilde J(s)\}.$  Set
$$ \ W(s):=\bigcup\limits_{j \in \widetilde J(s)}{\rm supp}(\mu(s,j)),\ \mu^*(s):=\frac{1}{|\widetilde J(s)|}\sum\limits_{j\in \widetilde J(s)}\mu(s,j)\ \forall s\in \widetilde S.$$
By construction,  we have:
$$\mu^*(s)\in \widetilde T_0(s)\subset T_0,\ {\rm supp}(\mu^*(s))=W(s)\ \forall s\in \widetilde S.$$
Since $\mu^*(s)\in  T_0,$ then it follows from (\ref{Ta}) that there exists $k(s)\in S$ such that
\be \mu^*(s)\in  T_0(k(s)) \ \Longrightarrow\ {\rm supp}(\mu^*(s))=W(s)\subset P_*(k(s)).\label{115}\ee
Thus, we have shown that
\be \forall s \in \widetilde S\ \ \exists \ k(s)\in S:\ W(s)\subset P_*(k(s)).\label{114}\ee
For $s\in S,$   consider  a vector $\tau^{*}(s)$ defined as follows:
$\tau^{*}(s):=\frac{1}{|J(s)|}\sum\limits_{j\in J(s)}\tau(j).$
By construction, $$\tau^{*}(s)\in T_0(s),\ \  {\rm supp}(\tau^{*}(s))=P_*(s),$$
 and, due to Proposition \ref{PP2}, it holds: $\tau^{*}(s)\in T_0$. Then it follows from  representation (\ref{z1})
that there exists $m(s)\in \widetilde S$ such that
\be \tau^{*}(s)\in \widetilde T_0(m(s)) \ \Longrightarrow\ {\rm supp}(\tau^{*}(s))=P_*(s)\subset  W(m(s))\ \forall s\in S.\label{15.1}\ee
Hence,  we  can  conclude that 
\be \forall s \in S \ \ \exists\, m(s)\in \widetilde S:\ P_*(s)\subset W(m(s)).\label{16}\ee

Suppose that there are two indices $s_1\in S$ and $s_2\in S$, $s_1\not=s_2,$ such that $m(s_1)=m(s_2)=m\in \widetilde S.$
Then, due to (\ref{16}) and Corollary \ref{C1}, we have
$$P_*(s_1)\subset W(m), \ P_*(s_2)\subset W(m), \ P_*(s_1)\not \subset P_*(s_2).$$
It follows from the  relations  above  and (\ref{114})   that
$$P_*(s_1)\cup P_*(s_2)\subset W(m)\subset P_*(k(m))\mbox{, where } P_*(s_1)\not\subset  P_*(s_2), \ k(m)\in S.$$
 This implies   that  in  the minimal zeros  graph $G(X)$, there exist  cliques $J(s_1)$ and $J(k(m))$ such that
$P_*(s_1)\subset P_*(k(m))$ and  $P_*(s_1)\not=  P_*(k(m))$. However,  this contradicts Corollary \ref{C1}.

Thus, we have shown that for any $s\in S$, there is an index $m(s)\in \widetilde S$ such that  the  inclusion
 in (\ref{16}) holds true and the
conditions $s \in S,$ $\bar s \in S,$ $s\not=\bar s,$ imply that $m(s)\not =m(\bar s).$
 These considerations permit us to conclude that $|\widetilde S|\geq |S|.$

Suppose now that $|\widetilde S|= |S|.$ It follows from  the reasoning   above that the mapping $m(s):S\to \widetilde S$ is an {\it one-to-one} mapping. Hence, without loss of generality, we may  consider that $\widetilde S= S$ and $m(s)=s$ for all $s \in S.$

  From  (\ref{114})  and  (\ref{16}), we have
$$P_*(s)\subset W(s)\subset P_*(k(s)) \ \mbox{ with some } k(s)\in S, \ \forall s \in S.$$
From the relations above  and Corollary \ref{C1}, it follows:
\be k(s)=s,\ P_*(s)=W(s) \ \forall s \in S.\label{*10}\ee
Then,
due to (\ref{115}) and (\ref{15.1}),  we obtain
\be \mu^*(s)\in T_0(s),\ \ \tau^{*}(s)\in \widetilde T_0(s)\ \forall s\in S.\label{200}\ee

Consider some $s\in S$ and suppose that there exists $i_0\in J(s)$ such that $\tau(i_0)\not \in \widetilde T_0(s).$
Since   $\tau(i_0)\in T_0(s)$, $\tau(i_0)\not \in \widetilde T_0(s)$, and $\mu^*(s)\in T_0(s)$,   it is evident that there exists $\alpha\in (0,1)$
such that
\be  \bar \tau:=\alpha\tau(i_0)+(1-\alpha)\mu^*(s)\in T_0(s),\ \bar \tau\not\in \widetilde T_0(s), \ {\rm supp}(\bar \tau)=P_*(s).\label{201}\ee
Here we took into account that ${\rm supp}(\mu^*(s))=W(s)=P_*(s)$ and ${\rm supp}( \tau(i_0))\subset P_*(s).$

On the other hand, since  $\bar \tau \in T_0$   and  $\bar \tau\not \in \widetilde T_0(s)$,  it follows from   (\ref{z1})  that there exists $\bar s\in \widetilde S=S$, $\bar s\not=s,$  such that
$\bar\tau\in \widetilde T_0(\bar s)$ and, hence, ${\rm supp}(\bar \tau)\subset W(\bar s).$ Taking into account
 this inclusion   and
relations (\ref{201}) and (\ref{*10}), we  obtain
$$ {\rm supp}(\bar \tau)=P_*(s)\subset W(\bar s)=P_*(\bar s) \mbox{ with } s \in S,\ \bar s\in S,\ s\not=\bar s.$$
 However, this contradicts (\ref{8.1}).

Thus we  have proved that
$\tau(i)\in \widetilde T_0(s)\ \forall i \in J(s),\ \forall s\in S$, and therefore
\be T_0(s)\subset \widetilde T_0(s) \ \forall s\in S.\label{4Z}\ee
In a similar way, one can show that $\mu(s,j)\in T_0(s)$ for all $j\in \widetilde J(s)$ and  all  $s \in S=\widetilde S.$
Hence, $\widetilde T_0(s)\subset  T_0(s)$ for all $s\in S.$ It follows from these inclusions and inclusions (\ref{4Z})  that
$\widetilde T_0(s)=T_0(s)$ for all $ s\in S.$ $\ \Box$

\vspace{3mm}

  Now we will give an alternative characterization/definition of the set $\{J(s), s\in S\}$ of subsets of $J$ producing the minimal representation of the
set $T_0.$

\begin{proposition}\label{PP1}    Consider a matrix  $X\in \cop$,   the corresponding  extended minimal zeros support set
  $ \mathcal{E}:=\{({\rm supp}(\tau(j)),M(j)), j \in J\}$,
and  the set (\ref{cc1})  of all maximal cliques of the corresponding  minimal zeros  graph $ G(X).$
Then the following conditions   are satisfied: \begin{enumerate}

\item[a)] $\ \ J=\bigcup\limits_{s \in S} J(s),$

\item[b)] $ \ \ P_*(s)\subset M(j)\ \forall j \in J(s),\ \forall s \in S,$

\item[c)] if $|S|\geq 2,$ then  for each pair of different indices $s\in S$ and $\bar s\in S$, it holds:
\be  J(s)\setminus J(\bar s)\not=\emptyset,\ J(\bar s)\setminus J( s)\not=\emptyset,  \label{c2.1}\ee
\be\forall i_0\in J(s)\setminus J(\bar s)\ \exists \ j_0\in J(\bar s)\setminus J( s) \mbox{ such that }\, {\rm supp}(\tau(i_0))\not \subset M(j_0).\label{c2.2}\ee
\end{enumerate}
\end{proposition}

{\bf Proof.} Let us  first prove   condition a). Suppose that on the contrary, there exists $j_0\in J$ such that $j_0\not \in \bigcup\limits_{s \in S} J(s).$
Consider a clique $\overline J=\{j_0\}$  and let $\overline {\overline J}$ be a maximal clique such that $\overline J\subset \overline {\overline J}$. Then, by construction,
$j_0\in \overline {\overline J}\setminus J(s)$ for all $ s \in S$. Consequently,
$\overline {\overline J}\not=J(s)$ for all $ s \in S.$ But this contradicts the
assumption that the set (\ref{cc1}) contains all maximal cliques of $G(X)$. Condition a) is proved.

 Now, let us prove  condition b). Suppose, on   the contrary, that  there exist $\bar s\in S$ and $j_0\in J(\bar s)$ such that
  $P_*(\bar s)\not\subset M(j_0)$. Hence, we can state that  
\be \exists \ k_0\in P_*(\bar s) \mbox{ such that } k_0\not \in M(j_0).\label{c6}\ee
By construction, $J(\bar s)$ is a clique    of the graph   $G(X)=(J,V)$.
 Therefore $(i,j)\in V$ for all $ i \in J(\bar s),$ $j \in J(\bar s),$ $ i<j,$ and, hence,
\be
 (\tau(i))^\top X\tau(j)=0\ \forall i \in J(\bar s),\ \forall j \in J(\bar s),\  i\leq j.\label{c7}\ee
  As  $ k_0\in P_*(\bar s)$,  it follows from (\ref{c2}) that there exists $i_0\in J(\bar s)$ such that $k_0\in  {\rm supp}(\tau(i_0)).$
Then, taking into account (\ref{c6}) and the definition of the set $M(j_0)$, we obtain
$$(\tau(i_0))^\top X\tau(j_0)=\sum\limits_{k\in  {\rm supp}(\tau(i_0))}\tau_k(i_0)e^\top_k X\tau(j_0)\geq \tau_{k_0}(i_0)e^\top_{k_0} X\tau(j_0)>0.$$
But this inequality contradicts (\ref{c7}).  Thus,  condition b) is proved.

 Finally, we will  prove  condition c). Suppose  that  there exist $s, \, \bar s\in S$  with  $s\not=\bar s$ such that
$J(s)\setminus J(\bar s)=\emptyset.$ This implies the inclusion  $J(s)\subset J(\bar s)$, which contradicts the assumption that $J(s)$ is a
 maximal clique  of  $G(X)$. Hence  conditions  (\ref{c2.1}) hold true.

To prove (\ref{c2.2}), suppose that there exists $i_0\in J(s)\setminus J(\bar s)$ such that
${\rm supp}(\tau(i_0))\subset M(j)\ \forall j \in J(\bar s)\setminus J( s).$
It follows from these inclusions and  condition b) proved above that
${\rm supp}(\tau(i_0)) \subset M(j)\ \forall j \in J(\bar s).$ This implies that for all $j \in J(\bar s)$, we have
$$\tau(i_0)^\top X\tau(j)=\sum\limits_{k\in {\rm supp}(\tau(i_0))}\tau_k(i_0)e^\top_k X\tau(j)=0.$$ Hence, by definition,
\be (i_0,j)\in V\ \ \forall j \in J(\bar s).\label{c8}\ee
Consider a set $\overline J(\bar s):=J(\bar s)\cup \{i_0\},$ $\overline J(\bar s)\subset J,$ $i_0\not \in  J(\bar s).$ It
follows from (\ref{c7}) and (\ref{c8}) that $\overline J(\bar s)$ is a clique    of  $G(X)$ which  contradicts the assumption that
$J(\bar s)$ is a maximal clique  of   $G(X)$. $\ \Box$

\begin{proposition}\label{PP4}
 Consider  $X\in \cop$ and the corresponding extended minimal zeros support set  $ \mathcal{E}= \{({\rm supp}(\tau(j)),M(j)), j \in J\}$.
Let  $\{J(s), s \in S\}$ be a set of  subsets of $J$ satisfying the
 conditions a)-c)  of  Proposition \ref{PP1} with   the  maximum  number $|S|$  of elements   and
the sets $ P_*(s), s \in S,$ defined in (\ref{c2}). Then the set $\{J(s), s \in S\}$ is the set of all maximal cliques of the graph $G(X).$
\end{proposition}
{\bf Proof.} It follows from  condition b) that for  $ s \in S,$ the set $J(s)$ is a clique  of  $G(X)$.
Let us show that for  $ s \in S,$ the clique  $J(s)$ is maximal.

Suppose the contrary: there exists $\bar s \in S$ such that the corresponding clique $J(\bar s)$ is not a maximal clique  of  $G(X)$.
 Hence, there exists $i_0\in J\setminus J(\bar s)$
such that $(i_0,j)\in V$ for all $j \in J(\bar s).$ Consequently, $\tau(i_0)^\top X\tau(j) =0$ for all $j \in J(\bar s).$
 From the   latter  equalities, we obtain the  inclusions
\be {\rm supp}(\tau(i_0))\subset M(j) \ \ \forall j \in J(\bar s).\label{*20}\ee
Note that it follows from condition a) that there is $s \in S$ such that $i_0\in J(s)$, and it is evident that $s\not=\bar s,$
$i_0\in J(s)\setminus J(\bar s)$ since $i_0\not \in J(\bar s).$
But the  conditions on the indices  $s$ and $\bar s$,   and   the  inclusions (\ref{*20}) contradict  the  condition c).
Thus we have shown that for any $s \in S,$ the  clique   $J(s)$ is a maximal  clique of  $G(X)$.

Now suppose that there is a maximal clique $J(s_0)$ of $G(X)$ such that $J(s_0)\not \in \{J(s), s \in S\}$,  $s_0\not\in S$.
 Let us denote $ \overline S:=S\cup\{s_0\}$
and consider the extended
set $\{J(s), s \in  \overline S\}$ of maximal cliques of $G(X)$. Then it is easy to show (see  the proof of Proposition \ref{PP1})
that this extended set satisfies  the  conditions a)-c)  which  contradicts the assumption that the set
$\{J(s), s \in S\}$ is maximum by  the  number $|S|$  of elements. $\ \Box$

\begin{corollary}\label{C4}  Consider  $X\in \cop$ and the corresponding extended minimal zeros support set
  $ \mathcal{E}= \{({\rm supp}(\tau(j)),M(j)), j \in J\}$.   Let  $\{J(s), s \in S\}$ be a maximal cardinality set of  subsets of $J$ satisfying  the
conditions a)-c)   in  Proposition \ref{PP1}  with the sets $ P_*(s), s \in S,$ defined in (\ref{c2}). Then $T_0=\bigcup\limits_{s\in S}T_0(s)$, where $T_0(s)={\rm conv}\{\tau(j), j \in J(s)\},$ $s \in S, $ is the minimal  representation of the set $T_0.$
\end{corollary}

{\bf Proof}.  The proof of this corollary  follows from Propositions \ref{PP2} and \ref{PP4}.   $\ \Box$

\vspace{2mm}

Corollary \ref{C4} generalizes  the  results  of     \cite{KT-Mathematics}, where it was shown that
if a set  $\{J(s), s \in S\}$ of  subsets of $J$ satisfies  the
conditions a)-c)   {of} Proposition \ref{PP1} with    the  maximum number $|S|$  of elements,  then  the set $T_0$ admits a  representation
 $T_0=\bigcup\limits_{s\in S}T_0(s)$ with $T_0(s)={\rm conv}\{\tau(j), j \in J(s)\},$
$s \in S. $  Notice that  it was not shown that this  {re}presentation is  minimal.
Also, in    \cite{KT-Mathematics}, there  were no explicit constructive rules provided to find the set    $\{J(s), s \in S\}$  of  subsets of $J$ satisfying  the
conditions a)-c)  {of} Proposition \ref{PP1}  and having the   maximum number $|S|$  of elements.
{Furthemore,} other properties of the set $\{J(s), s \in S\}$ were not characterized.

\vspace{2mm}

Due to   {the novel} results obtained in this  section, we  now  know that

\vspace{2mm}

$\bullet$   {the re}presentation of the set $T_0$   {described in \cite{KT-Mathematics}} is a {\it unique} minimal  {re}presentation{;}

$\bullet$   {the} set   $\{J(s), s \in S\}$ of  subsets of $J$ satisfying   conditions a)-c)  {of} Proposition \ref{PP1}  {and possessing}  the  maximum number $|S|$
 of elements is the set of maximal cliques  {in the minimal zeros} graph $G(X).$

\vspace{2mm}

The  findings of this paper  provide   new information about  the properties of the sets $J(s), s \in S,$
and,  importantly, offer explicit  rules for constructing these sets.

In fact, the problem of  determining  the list of all maximal cliques of a given undirected graph is    well known in graph theory.
There are many  algorithms  designed to construct the set of all maximal cliques  in a  given graph
 (see for example, \cite{1,2,3}). Hence, we can apply any one of these  algorithms to  construct
 the set $\{J(s), s \in S\}$ of  subsets of $J$ that generates the minimal {re}presentation (\ref{Ta}) of the set $T_0.$

 \vspace{2mm}

It follows from the results of this paper that
\begin{itemize} \item   the finite data set
$ \{\tau(j),\, j \in J;\ \ J(s), \, s \in S\}$
is the minimal data set that allows  one to completely describe  the set $T_0$ of all normalized zeros of $X\in \cop$;
 \item  there are algorithms that can be applied   to    construct this date set.
 \end{itemize}

\section{ Some remarks }
 In the previous sections,   it was  shown that  for a given matrix $X\in \cop$ with the corresponding  set $ \{\tau(j), j \in J \}$
 of  normalized   minimal  zeros,
 there exists a unique graph $G(X)$ such that the  corresponding  set $\{J(s), s \in S\}$ of all  its
 maximal cliques  defines the minimal  representation  of the set $T_0$ in the form $T_0=\bigcup\limits_{s\in S}T_0(s)$ with $T_0(s)={\rm conv}\{\tau(j), j \in J(s)\},$
$s \in S$.

 Notice that there may exist  several  copositive  matrices  having the same minimal zeros graphs, {\it e.g.} $X\in \cop$ and $Y\in {\cal COP}( p^*)$  such that $G(X)=G(Y)$  (see the  example  below). On the other hand, for any undirected graph  $G$,  there exists a copositive  matrix $Y$ such that $G=G(Y)$.  In fact, let us explicitly construct  such a matrix $Y$.


  Consider an  undirected graph  $G=(J, V)$ with the  set of vertices $J=\{1,\dots, p^*\}$  and  the set of edges $V$.
Consider a matrix $Y\in \mathbb R^{ p^*\times  p^*}$  with elements $Y_{ij}$, $i \in J,$ $ j \in J,$ defined by the rules
$$Y_{ii}=0 \ \forall i\in J;\ Y_{ij}=Y_{ji}=0 \mbox{ if } (i,j)\in V;\ Y_{ij}=Y_{ji}=1 \mbox{ otherwise}.$$
The matrix $Y$ is symmetric and all its elements are non-negative. Hence  $Y\in {\cal COP}(p^*)$.
Let us show that $G=G(Y)$, where $G(Y)=(J(Y), V(Y))$  is the minimal zeros graph for matrix $Y.$

By construction, for all $ i\in J,$ the vector $\bar\tau(i)=(\bar\tau_k(i), k\in  J)^\top$ with $\bar\tau_k(i)=0$ for $k \in J\setminus\{i\}$ and
$\bar\tau_i(i)=1$ is a normalized zero of $Y$. It is evident that all  vectors $\bar\tau(i), i \in J,$ are minimal normalized zeros of $Y$ and there does not exist
other minimal zero of $Y.$ Hence,   $\{\bar\tau(i), i \in J\}$ is the set of all normalized minimal zeros of $Y$.
Consequently,   $J(Y)=J.$

Now let us show that   the set $V(Y)$ of edges of the graph   $G(Y)$   coincides with  $V$:   $V(Y)=V.$
By construction, for $i \in J$ and $j\in J$, $i<j$, we have
$$ (\bar\tau(i))^\top Y\bar\tau(j)=Y_{ij}=\left\{\begin{array}{ll}
0 & \mbox{ if } (i,j)\in V,\cr
1 & \mbox{ otherwise. }\end{array}\right.$$
It follows from these relations that
$$V(Y):=\{(i,j): i\in J,\  j\in J,\ i<j, \ (\bar\tau(i))^\top Y\bar\tau(j)=0\}=V.$$

 For illustration, we will present below an example showing  that for different copositive matrices $X$ and $Y$, the equality $G(X)=G(Y)$ may take place.

\vspace{3mm}

\textbf{ Example.}
Consider matrices  $X$ and   $\overline{X}$ defined as follows:
{\small $$ X=\begin{pmatrix}
		0&0&1&1&1\cr
		0&0&1&1&1\cr
		1&1&0&0&1\cr
    1&1&0&0&1\cr
		1&1&1&1&1\end{pmatrix} \mbox{ and } \overline{X}=\begin{pmatrix}
1 &-1 &  1  & 1 &  -1\cr
    -1 &1  & -1 &  1 &   1.5 \cr
     1  &-1 & 1  & -0.5 &1.5\cr
     1 &  1 & -0.5 & 1 & -1\cr
    -1 & 1.5 & 1.5 &-1 & 1
		\end{pmatrix}=H+\Delta \overline{X}, $$
		$$ \mbox{where }	H=\begin{pmatrix}
1 &-1 &  1  & 1 &  -1\cr
    -1 &1  & -1 &  1 &   1 \cr
     1  &-1 & 1  & -1 &1\cr
     1 &  1 & -1 & 1 & -1\cr
    -1 & 1 & 1 &-1 & 1
		\end{pmatrix},
		\Delta \overline{X}=\begin{pmatrix}
0 &0 &  0  & 0 &  0\cr
    0 &0  & 0 &  0 &   0.5 \cr
     0  &0 & 0  & 0.5 &0.5\cr
     0 &  0 & 0.5 & 0 & 0\cr
    0 & 0.5 & 0.5 &0 & 0
		\end{pmatrix} .$$}

     For matrix $X$,  we have: $X\in {\cal COP}(5)\ $,
		$\ T_0=T_0(1)\cup  T_0(2),\ $ where \\ $\ T_0(1)={\rm conv}\{ \tau(1), \tau(2)\},\ T_0(2)={\rm conv}\{ \tau(3), \tau(4)\}$, and
		$$ \tau(1)=(1\, 0\, 0\, 0\,0)^\top,\  \tau(2)=(0\, 1\, 0\,  0\,0)^\top,\  \tau(3)=(0\, 0\,  1\, 0\,0)^\top,\
		 \tau(4)=( 0\, 0\, 0\, 1\,0)^\top$$
		are minimal zeros of $X$.
		 Hence, the sets forming the  corresponding extended minimal zeros support set  $ \mathcal{E}=\{({\rm supp}( \tau(j)),  M(j)), j =1,...,4\}$
		are as follows:
		$${\rm supp}(\tau(j))=\{i\} , \ i =1,\dots,4; \  M(1)= M(2)=\{1,2\},\   M(3)= M(4)=\{3,4\}.$$

\vspace{2mm}

		The  $5\times 5$ matrix  $H$  above   is called the {\it Horn matrix}  (see {\it e.g.} \cite{Hild2017,Hild2020,Vargas}, and the references therein) and it is
        known that  it is copositive.  It is easy to verify that $\Delta \overline{X}\in {\cal COP}(5)$ and hence
				we can conclude that $\overline{X}\in {\cal COP}(5)$.
		For  $\overline{X}$,   the corresponding set $\overline{T}_0$ of all normalized zeros takes the form
		$\overline{T}_0=\overline{T}_0(1)\cup  \overline{T}_0(2).$ Here $ \overline{T}_0(1)={\rm conv}\{\bar \tau(1),\bar \tau(2)\},\
		\overline{T}_0(2)={\rm conv}\{\bar\tau(3),\bar \tau(4)\}$,
			where
		$$\bar\tau(1)= {0.5}(1\, 1\, 0\, 0\, 0)^\top,\ \bar\tau(2)= {0.5}(0\, 1\, 1\, 0\, 0)^\top,\ \bar\tau(3)= {0.5}(0\, 0\, 0\, 1\, 1)^\top,\
		\bar\tau(4)= {0.5}(1\, 0\, 0\, 0\, 1)^\top$$
		are  normalized    minimal zeros of $\overline{X}$. The  following   sets form
		the  corresponding extended minimal zeros support set  $ \overline{\mathcal{E}}:=\{({\rm supp}( \bar \tau(j)),  \overline{M}(j)), j =1,...,4\}$:
		\bea &{\rm supp}(\bar \tau(1))=\{1,2\},\, {\rm supp}(\bar \tau(2))=\{2,3\},\, {\rm supp}(\bar \tau(3))=\{4,5\},\, {\rm supp}(\bar \tau(4))=\{1,5\},\nonumber\\
& \overline{M}(1)=\overline{M}(2)=\{1,2,3\},\, \overline{M}(3)=\overline{M}(4)=\{1,4,5\}.\nonumber\eea
		
		It is easy to see that $ G(X)=G(\overline{X})=(J,V)$ with $J=\{1,\dots,4\}$ and $V=\{(1,2),\, (3,4)\}.$ The set of maximal cliques of graph $G=(J,V)$ is
 	$\{J(s), s \in S\}$ with $S=\{1,2\},$  $J(1)=\{1,2\},$  and $J(2)=\{3,4\}.$

		Thus, we have shown that different copositive matrices  $X$ and $\overline{X}$ {generate} the same corresponding   minimal zeros  graph.
 Notice that  the  matrices  $X$ and $ \overline{X}$ have different
extended minimal zeros support sets
		${\mathcal{E}} $ and $\overline{\mathcal{E}}$, respectively.



\begin{thebibliography}{99}

\bibitem{Baumert66} Baumert L.D. (1966) Extreme copositive quadratic forms.   Pacific J. Mathematics. 19 (2): 197-204.

\bibitem{Bomze2012}  {Bomze I. M. (2012) Copositive optimization - Recent developments and applications. EJOR. 216(3): 509--520.}


\bibitem{1} Bron C. $\& $ Kerbosch, J. (1973) Algorithm 457: finding all cliques of an undirected graph. Communications of the ACM. 16(9): 575-577.

\bibitem{Dick2021} Dickinson P. J. $\&$ de Zeeuw R. (2021) Generating irreducible copositive matrices
using the stable set problem. Discrete Applied Mathematics, 296, 103-117.

\bibitem{Dick2016} Dickinson P. J. $\&$ Hildebrand R. (2016) Considering copositivity locally. Journal of Mathematical Analysis and Applications. 437(2): 1184-1195.

\bibitem{Dic2013} Dickinson P. J., D{\"u}r M., Gijben L. $\&$ Hildebrand R. (2013)
Irreducible elements of the copositive cone. Linear Algebra and its Applications. 439(6): 1605-1626.

\bibitem{Dic2019}   Dickinson P.J.C. (2019) A new certificate for copositivity. Linear Algebra
and its Applications. 569: 15-37.

\bibitem{Dur_2010}{D\"{u}r M. (2010) Copositive Programming -- a Survey. In: Diehl  M., Glineur F., Jarlebring E., Michielis W. (eds) Recent advances in Optimization and its applications in Engineering.  Springer-Verlag Berlin Heidelberg X1: 535 p.}


\bibitem{Hild2017} Hildebrand R. (2017) Copositive matrices with circulant zero support set. Linear Algebra and its Applications. 514: 1-46.

 \bibitem{H1} Hildebrand R. (2014) Minimal zeros of copositive matrices. Linear Algebra and its Applications. 459: 154-174.

 \bibitem{Hild2020} Hildebrand R. (2020) On the algebraic structure of the copositive cone. Optimization Letters. 14(8): 2007-2019.

\bibitem{2} Hua X., Zhong M., Liu Q. $\&$ Wang M. (2020) List all maximal cliques of an undirected graph:
 a parallable algorithm. In IOP Conference Series: Materials Science and Engineering. 790 (1): 012076. IOP Publishing.

\bibitem{Klerk} Klerk E. D. $\&$ Pasechnik D. V. (2002) Approximation of the stability number of a graph via
copositive programming. SIAM Journal on Optimization. 12(4): 875-892.


		\bibitem{KTD_2020} Kostyukova O.I., Tchemisova T.V. $\&$ Dudina O.S.  (2020)
			Immobile indices and CQ-free optimality criteria for linear copositive programming problems, Set-Valued Var. Anal. 28: 89-107.
			

\bibitem{KT_Opt_2022}    Kostyukova O.I.  $\&$    Tchemisova T.V.  (2022)
On equivalent representations and properties of faces of the cone of copositive matrices. Optimization. 71(11): 3211-3239.


\bibitem{KT-Mathematics}{Kostyukova O. $\&$  Tchemisova T. (2021)  Structural properties of faces of the cone of copositive matrices}.
Mathematics.   9(21): 2698. https://doi.org/10.3390/math9212698

\bibitem{3} Makino K. $ \&$ Uno T. (2004) New algorithms for enumerating all maximal cliques. In Algorithm Theory-SWAT 2004:
9th Scandinavian Workshop on Algorithm Theory, Humlebaek, Denmark. Proceedings. 9:  260-272.



\bibitem{Povh2007} Povh J. $\&$ Rendl F. (2007) A copositive programming approach to graph partitioning. SIAM Journal on Optimization. 18(1): 223-241.

\bibitem{Povh2013} Povh J. (2013) Contribution of copositive formulations to the graph partitioning problem. Optimization. 62(1): 71-83.

\bibitem{Vargas} Vargas L.F. $\&$  Laurent, M. (2023) Copositive Matrices, Sums of Squares and the Stability Number of a Graph.
In: Ko\v{c}vara M., Mourrain B., Riener C. (eds) Polynomial Optimization, Moments, and Applications.
 Springer, 206: 99-132.



\end{thebibliography}
\end{document}